\documentclass[11pt,reqno]{amsart}
\usepackage{amssymb,latexsym}
\usepackage{graphicx}
\usepackage{fancyhdr,amssymb}
\usepackage{url}		

\theoremstyle{plain}
\numberwithin{equation}{section}

\begin{document}
\fancyhead{}
\renewcommand{\headrulewidth}{0pt}
\fancyfoot{}
\fancyfoot[LE,RO]{\medskip \thepage}

\setcounter{page}{1}

\title[Galois Identities of the three term recurrence]{Galois Identities of the three term recurrence}
\author{Cheng Lien Lang}
\address{Department Applied of Mathematics\\
                I-Shou University\\
                Kaohsiung, Taiwan\\
                Republic of China}
\email{cllang@isu.edu.tw}
\thanks{}
\author{Mong Lung Lang}
\email{lang2to46@gmail.com}

\begin{abstract}
We study the identities associated to the characteristic polynomial $f(x)$ of certain
 recurrence relation. The existence of such identities is closely related to the
  fact that the Galois group of $f(x)$ over $\Bbb Q$ is $\Bbb Z_2$.

\end{abstract}

\maketitle

\vspace{-.8cm}

\section {Introduction}

\noindent Let $\{u_0, u_1, u_2, \cdots\}$ be a sequence defined by the following recurrence relation
$u_0 = 0, u_1= 1$,
 $$u_r= pu_{r-1} -qu_{r-2}. \eqno(1.1)$$

 \medskip
 \noindent
 where $p$ and $q$ are rational numbers. The characteristic polynomial of $
  \{u_r^n\}$ is defined to be (see for examples, [Br], [CK], [J], [K])

 $$\Phi_n(p, q, x)=
  \sum_{i=0}^{n}
  (-1)^{i} q^{i(i-1)/2} (n|i)_u x^i,\eqno(1.2)$$

 \medskip
 \noindent
 where $(n|k)_u$ is the generalised binomial coefficient (see Appendix A). Note that
$(n|0)_u=1$ and $(n|k)_u= u_nu_{n-1}\cdots u_{n-k+1}/u_ku_{k-1}\cdots u_1$ for $1\le k \le n$
 if $u_1 u_2 \cdots u_k \ne 0$.
  Let $\sigma$ and $\tau$ be roots of $x^2- px+q =0$. It is well known that

   $$\Phi_{n}(p,q, x) = \prod_{j=0}^{n} (x-\sigma^j\tau^{n-j}).\eqno(1.3)$$

 \medskip
 \noindent
 In this article,
 we shall associated to the above factorisation some identities, which we propose to call it
  the {\em Galois Identities} of $\Phi_n(p,q,x)$.
These identities  (see Proposition 3.4 and Corollary 3.6) must be well known among the experts. However,
our interpretation of the existence of these identities maybe  of some interest (see Discussion 3.7).

 \section {Fibonacci and Lucas numbers and $L_n^2 - 5F_n^2 = 4(-1)^n$}

\noindent
In (1.1), $p=1$ and $q=-1$ give the Fibonacci numbers.
 Applying (1.3), it is clear that
  the Galois group of $\Phi_n(1, -1, x) $ is $\Bbb Z_2$ and the
   splitting field of $\Phi_n(1,-1, x)$ is $\Bbb Q(\sqrt 5)$.
    Recall another well known fact about the factorisation of $\Phi_n(1,-1, x)$.

  $$\Phi_n(1,-1,x) = (-1)^n(x^2-L_nx +(-1)^n)\Phi_{n-2}(1,-1, -x),\eqno(2.1)$$

   \medskip
   \noindent where $L_n$ is the $n$-th Lucas number.
    Since $\Phi_n(1,-1,x)$ splits completely in $\Bbb Q(\sqrt 5)$ and $x^2-L_nx +(-1)^n =0$
   splits in $\Bbb Q(\sqrt{L_n^2 -4(-1)^n})$, one must have

   $$L_n^2 -4(-1)^n = 5A_n^2,\eqno(2.2)$$

\medskip
   \noindent for some $A_n\in \Bbb N$.
   This tells us that the difference between $4(-1)^n$ and the square of the Lucas number $L_n$ must be five times a square $A_n^2$. As for why $A_n$ must be $F_n$, the $n$-th Fibonacci number,
    we note that both $F_n^2$ and $L_n^2 $ satisfy the following recurrence relation

   $$X(n+3) = 2X(n+2)+2X(n+2)-X(n),\eqno(2.3)$$

\medskip
\noindent and that     $L_n^2-4(-1)^n$ and $F_n^2$ have the same initial values.
Hence we have just recover the following well known identity by investigating the splitting
 field of $\Phi_n(1,-1,x)$.

 $$L_n^2 -4(-1)^n = 5F_n^2.\eqno(2.4)$$

 \medskip
 \noindent {\bf Discussion 2.1.} The above investigation suggests that (2.4)
 is not just a numerical coincidence and
 can  be viewed as the consequence of the fact that  the Galois group of $\Phi_n(1,-1,x)$
  is $\Bbb Z_2$ and that the splitting field of $\Phi_n(1,-1,x)$ is
  $\Bbb Q(\sqrt 5) = \Bbb Q(\sqrt { L_n^2 + 4(-1)^n})$.
\medskip

\noindent {\bf Remark 2.2.} It follows easily from (2.4) that the following identities hold,
 the first identity is proved by Freitag and the second by Zeitlin and Filipponi (independently).
  See  [F] for more detail.

  $$\frac{\,L_n^2 - (-1)^a L_{n+a}^2\,}{F_n - (-1)^a F_{n+a}^2} = 5,\,\,
  \frac{\,L_n^2 + L_{n+2a}^2+ 8(-1)^n\,}{F_n + F_{n+2a}^2} = 5.\eqno(2.5)$$

\section {The Galois Identities of $w_r = pw_{r-1}-qw_{r-2}$}

\noindent In general, $\Phi_n(p,q, x)$ can be factorised completely into the following
     (Lemma 3.3 of Cooper and Kennedy [CK]).

    $$\Phi_{n}(p,q, x) = \prod_{j=0}^{n} (x-\sigma^j\tau^{n-j}),\eqno(3.1)$$

    \noindent where $\sigma $ and $\tau$ are roots of $x^2-px+q =0$. Hence
    the splitting field of $\Phi_n(p, q, x)$ is
     $\Bbb Q(\sqrt { p^2 - 4q})$
and
     the Galois group
     of $\Phi_n (p, q, x) $ is $  \Bbb Z_2$
     if and only if $p^2-4q$ is not a perfect square in $\Bbb Q$.
\medskip

\noindent {\bf Definition 3.1.} Let $\sigma$ and $\tau$ be given as
  in (1.3). Define $\{w_n\}$ to be the sequence  $w_n = \sigma^n +\tau^n$.

  \medskip
  \noindent {\bf Lemma 3.2.} {\em $w_0=2,$ $w_1= p$, $w_r= u_{r+1} -qu_{r-1}$ and
  $w_r= pw_{r-1}-qw_{r-2}$. Suppose that
   $p^2-4q\ne 0$. Then $u_r= (w_{r+1}-qw_{r-1})/(p^2-4q)$. }
   \medskip

   \noindent {\em Proof.} It is clear that $w_0=2$ and that $w_1=p$. Applying Binet's formula, one has
    $w_r= u_{r+1} -qu_{r-1}$ and that $w_r$  satisfies the recurrence
    $w_r= pw_{r-1} -qw_{r-2}$. Since

    $$w_r= u_{r+1} -qu_{r-1} = pu_r -2qu_{r-1}, \, w_{r-1} = u_r-qu_{r-2}= 2u_r -pu_{r-1},\eqno(3.2)$$

 \medskip
 \noindent we conclude that   $u_r= (w_{r+1}-qw_{r-1})/(p^2-4q)$. This completes the proof of the lemma. \qed

\medskip
\noindent {\bf Lemma 3.3.} {\em  $q^n$,  $w_n^2$ and $u_n^2$ satisfy the following recurrence

 $$X(m+3) = (p^2-q)X(m+2) +(q^2-p^2q)X(m+1) +q^3X(m).\eqno(3.3)$$}

\noindent
{\em Proof.} Applying (1.1), one can show easily that
 $u_n^2$ satisfies the recurrence (3.3). The rest can be verified
  similarly. \qed

\medskip
\noindent
Suppose that $p^2 - 4q$ is not a square in $\Bbb Q$.
 Applying Galois Theory, the set of conjugates of $\sigma^n$ over $\Bbb Q$ is $\{\sigma^n, \tau^n\}$.
 Hence the following holds for every $n \in \Bbb N$.
 $$ f_n(x) = (x-\sigma^n)(x-\tau^n)
 = x^2 - w_n x + q^n
 \in \Bbb Q[x].\eqno (3.4)$$

\medskip
\noindent
It is clear that $f_n(x)$ splits in $\Bbb Q(\sqrt { w_n^2-4q^n})$.
 Applying (3.1),  $f_n(x)$ splits in  $\Bbb Q(\sqrt {p^2-4q})$.
  Hence
  $$ w_n^2 - 4q^n = z^2(p^2-4q),\eqno(3.5)$$

  \medskip
  \noindent for some $z \in \Bbb Q$. The following proposition shows that the solution of (3.5) is $z = u_n$,
   the recurrence we defined in (1.1).

\medskip
\noindent {\bf Proposition 3.4.} {\em Let $ w_n $ be given as in Lemma $3.1$. Then $ w_n^2 - 4q^n = u_n^2(p^2-4q).$}

\medskip
\noindent {\em Proof.}
Since both $ w_n^2 - 4q^n$   and $u_n^2(p^2-4q)$ satisfy the recurrence (3.3)
 and admit the same initial values, we have
$ w_n^2 - 4q^n = u_n^2(p^2-4q).$ \qed
\medskip

\noindent {\bf Corollary 3.5.} {\em
The equation $x^2 +y^2 -z^2 = 4$ is solvable in $\Bbb Z$.
Further, one may choose $y$ and $z$ in such a way that $py = 2z$ for any $p \in \Bbb Z$.
 }

\medskip

\noindent {\em Proof.} Let $q=1$ and let $p$ be any integer. Applying  Proposition 3.4, one has
$w_n^2 +(2u_n)^2 - (pu_n)^2 = 4$ for all $n \ge 1$.\qed

\medskip
\noindent {\bf Corollary  3.6.} {\em$ w_{2n}- 2q^n = u_n^2(p^2-4q).$
In particular, $L_{2n} -2(-1)^n = 5F_n^2$.
}

\medskip
\noindent {\em Proof.} Applying Definition 3.1, one has $w_n^2 = w_{2n} +2q^n$. \qed

\medskip
\noindent {\bf Discussion 3.7.}
The identities  in Proposition 3.4  and Corollary 3.6 must be  well known and we propose to call it the {\em Galois identities} associated to
 $\Phi_n(p, q, x)$. We would like to emphasise that  Proposition 3.4 is not just a numerical
  coincidence but can  be treated as the consequence of the fact
   that
the splitting field of $\Phi_n(p,q,x)$ is
  $\Bbb Q(\sqrt {p^2-4q}) = \Bbb Q(\sqrt { w_n^2 -4q^n})$. The most famous identity among all, of
   course, is (2.4).

\section {Appendix A}
  Let $\sigma$ and $ \tau$ be roots of $x^2 -px+q=0$ (known as the characteristic polynomial
  of (1.1)). It is well known that
 $$u_r = \sum_{i=1}^{r-1} \sigma^{r-1-i}\tau^i = f_r(\sigma, \tau) ,\eqno(A1)$$

 \medskip
\noindent where $f_r(x,y)$ is the polynomial  $ \sum x^{r-1-i}y^i  =(x^r-y^r)/(x-y).$
 Let $G(z)$ be the following function.
 $$ G(z) = (1-z^{m+n})(1-z^{m+n-1})\cdots (1-z^{m+1}/(1-z^n)(1-z^{n-1})\cdots (1-z).\eqno(A2)$$
 $G(z)$ is known as the Gaussian binomial coefficient  and is a polynomial in $z$ (a more powerful result actually implies that (A2) can be written
  as product of cyclotomic polynomials).
  One may apply this fact to show that
  that
 $(f_r f_{r-1}\cdots f_{r-k+1})/(f_kf_{k-1}\cdots f_1)\in \Bbb Z[x,y]$ is  a polynomial in $x$ and $y$.
Denoted by $F(r,k, x, y)$ this polynomial. Define

 $$ (r|k)_u = F(r,k,\sigma, \tau).\eqno(A3)$$

\medskip
\noindent We call $(A3)$ the generalised binomial coefficient. It is clear that if $u_1 u_2 \cdots u_k
 \ne 0$. Then $(n|k)_u$ takes a better looking form.

 $$(n|k)_u= u_nu_{n-1}\cdots u_{n-k+1}/u_ku_{k-1}\cdots u_1. \eqno(A4)$$

\bigskip

\noindent MSC2010 : 11B39, 11B83.

\medskip
\noindent fibonacci-7-2.tex

\end{document}